\documentclass[12pt,twoside,a4paper]{article}
\usepackage{graphicx}
\usepackage{amsfonts}
\usepackage{bbm}
\usepackage{mathrsfs}
\usepackage{fancyhdr}
 \usepackage{mathrsfs,amsmath,amssymb,amsthm}
 \usepackage{amssymb}
 \usepackage{amsbsy}
 \usepackage{color}
\allowdisplaybreaks

 \theoremstyle{plain}
\theoremstyle{remark}  \newtheorem{remark}{\noindent\mbox{Remark}}
 \theoremstyle{plain}
 \theoremstyle{plain}\newtheorem{lemma}{\noindent\mbox{Lemma}}
\theoremstyle{plain} \newtheorem{theorem}{\noindent\mbox{Theorem}}
 \theoremstyle{plain}\newtheorem{proposition}{\noindent\mbox{Proposition}}
 \theoremstyle{plain}\newtheorem{corollary}{\noindent\mbox{Corollary}}
\theoremstyle{definition} 

 \def\qed{\hfill$\Box$\medskip}
 \def\rto{\rightarrow\infty}
\def\z{\left}
\def\y{\right}
 \def\no{\nonumber}

\begin{document}
 \title{{On a maximum of nearest-neighbor random walk with asymptotically zero drift on lattice of positive half line}\footnote{Supported by National
Nature Science Foundation of China (Grant No. 11501008; 11601494).}}                  

\author{Hongyan \uppercase{Sun}\footnote{Email:sun$\_$hy@cugb.edu.cn; Department of Mathematics, China University of Geosciences, Beijing, 100083, China}~ $\ \&\ $
 Hua-Ming \uppercase{Wang}\footnote{Email:hming@ahnu.edu.cn; Department of Statistics, Anhui Normal University, Wuhu, 241003, China }}
\date{}
\maketitle%

\vspace{-.5cm}

\begin{center}
\begin{minipage}[c]{12cm}
\begin{center}\textbf{Abstract}\quad \end{center}
Consider a nearest-neighbor random walk with certain asymptotically zero drift on the positive half line. Let $M$ be the maximum of an excursion starting from $1$ and ending at $0.$ We study the distribution of $M$ and characterize its asymptotics, which is quite different from those of simple random walks.
\vspace{0.2cm}


\end{minipage}
\end{center}

\section{Introduction}\label{intro}
Let $X=\{X_k\}_{k\ge 0}$ be a Markov chain on $\mathbb Z^+:=\{0,1,2,...\}$ starting from some $x_0\in\mathbb Z_+$ and with transition probabilities
\begin{align}
  &P(X_{k+1}=1|X_k=0)=p_0=1\no\\
  &P(X_{k+1}=n+1|X_k=n)=p_n\in(0,1),\no\\
 &P(X_{k+1}=n-1|X_k=n)=q_n:=1-p_n,n\ge 1.\no
\end{align}
Unless otherwise stated, we always assume $x_0=1.$
For $i\ge 1$ set $$\rho_{i}=\frac{q_i}{p_i}{ ,}$$ which we will work with throughout.

Next we define the so-called maximum we concern for the chain $X.$  Let
\begin{align}
  &D=\inf\{k\ge 1: X_k<X_0\},\no\\
  &M=\sup\{X_k:0\le k\le D\},\no
\end{align}
where and throughout, we use the convention $\inf \phi=\infty.$ Clearly, $D$ is
the time that the walk $X$ hits   $X_0-1$ for the fist time and $M$ is the maximum that the walk  has ever reached in the time interval $[0,D).$

Our aim is to compute the exact distribution { of} $M$ and study its asymptotics. If for $n\ge1,$
$p_n\equiv p\in(0,1),$ then $X$ is just a simple random walk and in this case for all $i\ge1,$  $\rho_i\equiv \rho:=q/p.$ In this situation, the   asymptotical  distribution of $M$ on the event $\{D<\infty\}$ is as the following proposition.
\begin{proposition}\label{sm}
  Suppose that $p_i\equiv p\in (0,1), i\ge1$ and let $\rho:=\frac{1-p}{p}.$ Then  \begin{align}\label{asms}
  P(M=n,D<\infty)\sim \left\{\begin{array}{ll}
                              \frac{1}{n(n+1)}, & \rho=1, \\
                               (1-\rho)^2\rho^n, & \rho<1, \\
                               (1-\rho)^2\rho^{-(n+1)}, & \rho>1 ,
                             \end{array} \right.
 \text{ as } n\rto.
\end{align}
\end{proposition}
\noindent Proposition \ref{sm} has nothing new and can be easily  deduced from the following Corollary \ref{cor1}. Clearly, if $X$ is recurrent, then $P(D<\infty)=1.$
Taking the well-known recurrence criterion for simple random walk into account, we come to the conclusion that, as $n\rto,$

 \noindent (a) if $\rho=1,$  $X$ is {\it null recurrent} and $P(M=n)$ decays  {\it polynomially};

 \noindent (b) if $\rho<1,$  $X$ is {\it transient} and $P(M=n,D<\infty)$ decays {\it exponentially};

 \noindent  (c)  if $\rho>1,$  $X$ is {\it positive recurrent} and $P(M=n)$ decays {\it exponentially}.

One asks naturally that, besides the above polynomial and exponential speeds,   are there any situations, under which  $P(M=n,D<\infty)$ decays with other speed, for example, $n^{-s},s> 0$  or $ (n\log^\beta n)^{-1},\beta>1,$ et al?

The answer is affirmative. Adding some perturbation on the recurrent random walk, we get a near-recurrent random walk,  known  as Lamperti random walk which dates back to Harris \cite{har} and Lamperti \cite{lam} and is  extensively studied by many other authors, refer for example, to \cite{cfr,cfrb,hmw13,jlp,lo,w19} etc. To introduce Lamperti random walk,
we borrow the perturbations from \cite{cfr}.

For $K\ge 1,\ B\in R,$ set
\begin{align*}
    &\Lambda(1,i,B)=\frac{B}{i},\\
   &\Lambda(2,i,B)=\frac{1}{i}+\frac{B}{i\log i},\dots,\\
    &\Lambda(K,i,B)=\frac{1}{i}+\frac{1}{i\log i}+...+\frac{1}{i\log i\cdots\log_{K-2}i}+\frac{B}{i\log i\cdots\log_{K-1}i},
  \end{align*}
  where $\log_0 i=i, \log_1 i=\log i,.., \log_{K}i=\log \log_{K-1} i.$
 For $K$ and $B$ fixed, set $$i_0=\min\z\{i: \log_{K-1}i>0, \frac{|\Lambda(K,i,B)|}{4}< \frac{1}{2} \y\}.$$
 For fixed $B\in \mathbb R$ and $K=1,2,...$ set
\begin{align*}r_i=\left\{\begin{array}{ll}
   \frac{\Lambda(K,i,B)}{4}, & i\ge  i_0, \\
 r_{i_0}, & i< i_0.
\end{array}\right.\end{align*}


The following fact, parts of which can be found in \cite{cfr}, page 208, is a corollary of Lemma \ref{prhe} and Lemma \ref{ad} below.

\noindent{\bf Recurrence criterion.} {\it For $K\ge 2,$ when $p_i=\frac{1}{2}+r_i,\ i\ge 1,$
  the chain $X$ is recurrent if $B\le 1$ and transient if $B>1;$  when $p_i=\frac{1}{2}-r_i,\ i\ge 1,$ the chain $X$ is positive recurrent if $B>1$ and null recurrent if $B\le 1.$
For $K=1,$ when $p_i=\frac{1}{2}+r_i,\ i\ge 1,$
  the chain $X$  is transient, recurrent or positively recurrent according  as $B>1,$ $-1\le B\le 1,$ or $B< -1$ respectively.}

In the rest of the paper, unless otherwise stated, $ c$  is a strictly positive constant which may change from line to line.

\begin{theorem}\label{mxt}   Fix $K\ge1$ and $B\in \mathbb R.$
 {\rm (i)} If
    $p_i=\frac{1}{2}+r_i,i\ge 1,$
then, as $n\rto,$
\begin{align}\label{ap}
  P(M=n,D<\infty)\sim \left\{ \begin{array}{ll}
                                \frac{c }{n\log n\log\log n\cdots \log_{K-2}n\log_{K-1} n(\log_K n)^2}, &  \text{if } B=1,  \\
\frac{c}{n\log n\log\log n\cdots \log_{K-2}n(\log_{K-1} n)^B}, &  \text{if } B>1,  \\
\frac{c}{n\log n\log\log n\cdots \log_{K-2}n(\log_{K-1} n)^{2-B}}, &  \text{if } B<1.
                              \end{array}
  \right.
\end{align}
{\rm(ii)} If
$  p_i=\frac{1}{2}-r_i,i\ge 1,$
then, as $n\rto,$
\begin{align}\label{an}
  P(M=n,D<\infty)\sim
  \left\{\begin{array}{lll}
&\frac{c}{n^{B+2}}, &\text{if } K=1,B>-1,\\
&\frac{c}{n(\log n)^2}, &\text{if } K=1, B=-1,\\
 & cn^B, &\text{if } K=1,B<-1,\\
 &\frac{c}{n^3\log n...\log_{K-2}n (\log_{K-1}^n)^B}, &\text{if } K>1.\\
\end{array}\right.
\end{align}

\end{theorem}
\begin{remark} I) If $K=1,$ the result in (\ref{ap}) should overlap the one in (\ref{an}).
 Indeed, in this case,  since $\Lambda(1,i,B)=\frac{B}{i},$ replacing $B$ by $-B$ in (\ref{ap}), we get \eqref{an} and vice versa. II) We see that the asymptotics of $M$ are quite different from those of the simple random walk in (\ref{asms}). For example, fix $K=1$  and let $p_i=1/2+r_i,i\ge 1.$ Let us consider the null recurrent case. That is $-1\le B\le 1.$
  If $B=1,$  $P(M=n)$ decays with speed $\frac{c}{n\log^2 n},$   but if $-1\le B<1,$  $P(M=n)$ decays with speed $\frac{c}{n^{2-B}}.$   We conclude that  even for null recurrent case, the decay speeds are quite  different among different $B.$    While,  for simple null-recurrent random walk, $P(M=n)$ always decays with polynomial speed $cn^{-2}.$
\end{remark}
The rest of the paper is arranged as follows. In Section \ref{pr}, firstly,  we describe the distribution of $M$ by the escape probability, which can be written as the function of $\rho_1\cdots\rho_j, j\ge 1.$ Secondly, we study the limit behaviors of $\rho_1\cdots\rho_n$ as $n\rto$, which are essential for proving Theorem \ref{mxt} and the recurrence criterion. Finally, Theroem \ref{mxt} is proved in Section \ref{p3}.
\section{Preliminaries}\label{pr}

For $0\le a\le k\le b,$ let
\begin{align*}
  P_k(a,b,-)=P(X\text{ hits }a \text{ before\ } b |X_0=k).
\end{align*}
It is easily seen that \begin{align}\label{inv}
  P_a(a,b,-)=1,\  P_b(a,b,-)=0.
\end{align}
By Markov property, we have for $a<k<b,$
\begin{align}\label{hm}
  P_k(a,b,-)=p_kP_{k+1}(a,b,-)+q_kP_{k-1}(a,b,-).
\end{align}
Solving the homonic system (\ref{hm}) with the boundary condition (\ref{inv}), we get the following lemma, which is very standard and can found in many   documents, see for example \cite{sol,ze}.
\begin{lemma} \label{esx}For $0\le a\le k\le b,$ we have
\begin{align}P_k(a,b,-)=\frac{\sum_{j=k}^{b-1}\rho_{a+1}\cdots \rho_{j}}{1+\sum_{j=a+1}^{b-1}\rho_{a+1}\cdots \rho_{j}}.\label{ep}\end{align}
  \end{lemma}
\noindent Here and   throughout  the paper, we always assume that empty product equals $1$ and empty sum equals 0.

On the event $\{D<\infty\},$ the distribution of $M$ can be deduced from  Lemma \ref{esx} directly.
\begin{corollary}\label{cor1} For the chain $X,$   we have
\begin{align*}
  P(M=n,&D<\infty)=(1-P_1(0,n,-))P_n(0,n+1,-)\\
  &=\frac{1}{1+\sum_{j=1}^{n-1}\rho_{1}\cdots \rho_{j}}\cdot\frac{\rho_1\cdots \rho_n}{1+\sum_{j=1}^{n}\rho_{1}\cdots \rho_{j}}.
\end{align*}
\end{corollary}

Both the asymptotics of $M$  and the recurrent criterion of $X$ rely on the limit behavior of $\rho_1\cdots \rho_n,$ as $n\rto,$ which we   will give in the next lemma.
\begin{lemma}\label{prhe} Fix $K=1,2,..$ and $B\in \mathbb R.$ {\rm(a)} If $p_i=\frac{1}{2}+r_i, i\ge 1,$ then
\begin{align}\label{ra}
  \rho_1\cdots\rho_n\sim  \frac{c}{n\log n\log\log n\cdots \log_{K-2}n(\log_{K-1} n)^B}, \mbox{as}\ n\rto.
\end{align}
{\rm(b)} If $p_i=\frac{1}{2}-r_i,i\ge1,$ then
\begin{align*}
  \rho_1\cdots\rho_n\sim  c{n\log n\log\log n\cdots \log_{K-2}n(\log_{K-1} n)^B}, \mbox{as}\ n\rto.
\end{align*}
\end{lemma}

Before proving Lemma \ref{prhe}, let us study the recurrent criterion of $X$ at first.
Define  another Markov chain $X'=\{X'_k\}_{k\ge0},$ which starts from some $x_0'\in \mathbb Z_+,$ with transition probabilities
\begin{align}
  &P(X'_{k+1}=1|X'_k=0)=p'_0=1\no\\
  &P(X'_{k+1}=n+1|X'_k=n)=p'_n\in(0,1),\no\\
 &P(X_{k+1}=n-1|X_k=n)=q'_n:=1-p'_n,n\ge 1.\no
\end{align}
In literatures, if $q'_n\equiv p_n$ for all $n\ge 1,$ $X'$ is called the adjoint chain of $X$ and vise versa. We have the following recurrence criterion, which can  be found in Derriennic \cite{dey98}, page 204-205.
\begin{lemma}\label{ad}
  {\rm i)} The chain $X$ is positive recurrent if and only if its adjoint chain $X'$ is transient and vice versa.
     {\rm ii)} Both adjoint chains $X$ and $X'$ are null recurrent simultaneously.
\end{lemma}

It is easily seen from (\ref{ep}) that $X$ is transient or recurrent according as $\sum_{j=1}^{\infty}\rho_1\cdots \rho_j<\infty \text{ or }=\infty$ respectively.
Thus, with Lemma \ref{ad} in hands, the recurrence criterion for $X$ in Section \ref{intro}  follows from  (\ref{ra}) since $\sum_{j=1}^{\infty}\rho_1\cdots \rho_j$ converges if and only if $B>1.$

\vspace{0.5cm}

\noindent{\it Proof of Lemma \ref{prhe}.}
{ Case 1: Assume $p_i=1/2+r_i, i\ge 1.$} Then,
\begin{eqnarray*}
\rho_n=\frac{q_n}{p_n}=\frac{\frac12-r_n}{\frac12+r_n}=1-4r_n+8r^2_n+o(r^2_n), \ \mbox{as}\ n\rightarrow \infty.
\end{eqnarray*}
By the mean value theorem,  there exists    $\theta_n$ between 0 and $\rho_n-1$
such that
$$\log(\rho_n)=(\rho_n-1)-\frac1{2(1+\theta_n)^2}(\rho_n-1)^2.$$
Since $0<\frac1{2(1+\theta_n)^2}(\rho_n-1)^2<\frac{(\rho_n-1)^2}2\sim\frac{1}{2n^2}$ as $n\rightarrow\infty,$  $\sum\limits_{n=1}^{\infty}\frac1{2(1+\theta_n)^2}(\rho_n-1)^2$ is convergent.
Thus, \begin{align}\label{lim1}
\lim\limits_{n\rightarrow\infty}
\Big[\sum_{i=1}^n\log\rho_i-\sum_{i=1}^n(\rho_i-1)\Big]\mbox{\ exists.}
\end{align}
Since $\rho_n-1+4r_n=8r^2_n+o(r^2_n)$ and  $r^2_n\sim\frac1{16n^2}$ as $n\rto,$
 \begin{align}\label{lim2}
\lim\limits_{n\rightarrow\infty}\Big[\sum_{i=1}^n(\rho_i-1) +\sum_{i=1}^n(4r_i) \Big]\mbox{\ exists.}
\end{align}

Let $f(x):=\frac1x+\frac1{x\log x}+...+\frac{1}{x\log x\cdots\log_{K-2}x}+\frac{B}{x\log x\cdots\log_{K-1}x}$. Then $f(n)=4r_n, \forall n\ge {i_0}$, and for $K>1$, or $K=1,\ B\ge 0$, there exists an integer, say $M(M\ge i_0)$, which only depends on $B$ and $K$, such that  $f(x)$ is decreasing in $[M,\infty)$. Therefore, we have
$$\sum _{k={M}+1}^nf(k)\le \int_{M}^n f(x)dx\le\sum _{k={M}}^{n-1}f(k),\ \forall \ n>{M},$$
and then  \begin{align}\label{bd}
 0\le f(n)\le \sum _{k={M}}^nf(k)-\int_{M}^n f(x)dx\le f({M}),\ \forall \ n>{M}.
\end{align}
Since $\sum_{k={M}}^{n+1}f(k)-\int_{M}^{n+1} f(x)dx-[\sum_{k={M}}^nf(k)-\int_{M}^n f(x)dx]=f(n+1)-\int_n^{n+1}f(x)dx\le 0,\ \forall \ n>{M}.$  Then the sequence $\{\sum_{k={M}}^nf(k)-\int_{M}^n f(x)dx\}_{n\ge {M}},$ is decreasing. Therefore, owing to (\ref{bd}), we have that $\sum_{k={M}}^nf(k)-\int_{M}^n f(x)dx$ converges to some constant as $n\rto$. That is,
\begin{eqnarray}\label{lim3} \lim\limits_{n\rightarrow\infty}\Big[\sum_{i={M}}^n(4r_i)- \int_{M}^n \frac1x+...+\frac{B}{x\log x\cdots\log_{K-1}x}dx\Big]\mbox{\ exists}.
\end{eqnarray}
When $K=1, B<0$, let $f(x):=-\frac{B}{x}$, a similar argument as above  also  yields \eqref{lim3}.
From \eqref{lim1}, \eqref{lim2} and \eqref{lim3}, we get
 \begin{eqnarray*}
&&\frac{\rho_1...\rho_n}{\exp\{-\int_{M}^n \frac1x+...+\frac{B}{x\log x\cdots\log_{K-1}x}dx\}}\\
&&\quad =\exp\Big\{\sum_{i=1}^n(\log\rho_i)+\int_{M}^n \frac1x+...+\frac{B}{x\log x\cdots\log_{K-1}x}dx\Big\}
\end{eqnarray*}
converges to some constant as $n\rto$.
But
\begin{align*}
&\int_{M}^{n}\Big(\frac1x+\frac1{x\log x}+...+\frac{1}{x\log x\cdots\log_{K-2}x}+\frac{B}{x\log x\cdots\log_{K-1}x}\Big)dx\\
&=\log n+\log\log n+...+\log_{K-1}n +B\log_K n\\
&\qquad\qquad-(\log {M}+\log\log {M}+...+\log_{K-1}{M} +B\log_K {M}).
\end{align*}
Therefore,
\begin{eqnarray*}
\rho_1...\rho_n\sim\frac c{n\log n...\log_{K-2}n(\log_{K-1}n)^B}, \ \mbox{as} \ n\ \rto.
\end{eqnarray*}
Case 2. Assume $p_i=1/2-r_i, i\ge 1.$ The result can be deduced directly from Case 1. 
\qed

\section{Proof of Theorem \ref{mxt}}\label{p3}

Case 1. Assume $p_i=1/2+r_i, i\ge 1.$  When $B>1$,
\begin{align*}
\int_{i_0}^n&\frac1{x\log x...\log_{K-2}x(\log_{K-1}x)^B}dx=\int_{i_0}^n\frac1{(\log_{K-1}x)^B}d\log_{K-1}x\\
&\quad=\frac{1}{(B-1)}(\log_{K-1}{i_0})^{1-B}-\frac{1}{B-1}(\log_{K-1}n)^{1-B}\\
&\quad\rightarrow \frac{1}{(B-1)}(\log_{K-1}{i_0})^{1-B},\ \mbox{as} \ n\ \rto.
\end{align*}
So $\sum\limits_{n={i_0}}^{\infty}\frac1{n\log n...(\log_{K-1}n)^B}$ is convergent.
Thus, by Lemma \ref{prhe}(a), we get that $\sum\limits_{n=1}^{\infty}\rho_1...\rho_n$ is convergent.
By Corollary \ref{cor1}, we have
\begin{align*}
P(M&=n,D<\infty)
=\frac1{1+\sum_{j=1}^{n-1}\rho_1...\rho_j}\cdot\frac{ \rho_1...\rho_n}{1+\sum_{j=1}^{n}\rho_1...\rho_j},\no\\
&\sim  c\rho_1...\rho_n\sim  \frac c{n\log n...\log_{K-2}n(\log_{K-1}n)^B}\text{ as } n\rto.
\end{align*}
When $B\le1,$ for any $x_0\ge i_0,$
\begin{align}\label{int1}
& \int_{x_0}^n\frac1{x\log x...\log_{K-2}x(\log_{K-1}x)^B}dx=\int_{x_0}^n\frac1{(\log_{K-1}x)^B}d\log_{K-1}x \nonumber\\
&\quad\quad=\left\{
\begin{array}{rcl}
&\frac{1}{1-B}(\log_{K-1}n)^{1-B}-\frac{1}{(1-B)}(\log_{K-1}x_0)^{1-B},& B<1,\no \\
&\log_Kn-\log_Kx_0, & B=1,
\end{array} \right.\\
&\quad\quad\rightarrow\infty \text{ as } n\rto.
\end{align}
Note that when $K>1$ or $K=1,B\ge 0,$  $x\mapsto\frac1{x\log x...\log_{K-2}x(\log_{K-1}x)^B}$ is a decreasing function in $[n_0, \infty),$ for some integer $n_0\ge i_0$ depending on $K\text{ and }B.$
Therefore,  by a similar argument as the proof of \eqref{lim3}, we get
\begin{align}\label{int2}
\sum\limits_{i=n_0}^{n}\frac1{i\log i...(\log_{K-1}i)^B}\sim\int_{n_0}^n\frac1{x\log x...(\log_{K-1}x)^B}dx,\text{ as }n\rto.
\end{align}
When $K=1,B<0$,
Let $g(x):= \frac1{x^B}$. Then $g(x)$ is increasing in $[i_0,\infty)$.
We have $$g(n_0)\le \sum\limits_{i=i_0}^n g(k)-\int_{i_0}^n g(x)dx\le g(n).$$
Note that $g(n)=o(\int_{i_0}^n g(x)dx)$ as $n\rto$. Hence,
\begin{align}\label{int3}
\sum\limits_{i=i_0}^n \frac1{i^B}\sim\int_{i_0}^n  \frac1{x^B} dx,\text{ as }n\rto,
\end{align}
 which coincides  with \eqref{int2} whenever $K=1, B<0$.

From Corollary \ref{cor1}, \eqref{int1}, \eqref{int2} and \eqref{int3}, it follows that for $K\ge1$ and  $ B\le 1$,
\begin{align*}
P(&M=n,D<\infty)
=\frac1{1+\sum_{j=1}^{n-1}\rho_1...\rho_j}\cdot\frac{ \rho_1...\rho_n}{1+\sum_{j=1}^{n}\rho_1...\rho_j}\\
&\sim\frac{c\rho_1...\rho_n}{\z(\sum_{j=1}^{n}\rho_1...\rho_j\y)^2}\sim \left\{
\begin{array}{rcl}
&\frac{c}{n\log n...\log_{K-2}n(\log_{K-1}n)^{2-B}},& B<1\\
&\frac{c}{n\log n...\log_{K-1}n(\log_Kn)^2}, & B=1,
\end{array} \right.
\end{align*}
as $n\rto.$

Case 2. Assume $p_i=\frac{1}{2}-r_i,i\ge1.$  Set $e_2^x:=e^{e^x}$  and let $e_{k+1}^x:=e^{e_{k}^x},\ \forall k\ge2$.

 When $K>1,$ changing the variable in the integral, we have
\begin{eqnarray*}
\int_{i_0}^n {x\log x...(\log_{K-1}x)^B}dx&=&\int_{\log_{K-1}i_0}^{\log_{K-1}n}(e_{K-1}^ye_{K-2}^y...e^y)^2y^Bdy.
\end{eqnarray*}
On the other hand,
\begin{align*}
&\lim\limits_{x\rto}\frac{ \int_{\log_{K-1}i_0}^{\log_{K-1}x}(e_{K-1}^ye_{K-2}^y...e^y)^2y^B dy}{x^2\log x...\log_{K-2}x(\log_{K-1}^x)^B}\\
&=\lim\limits_{x\rto}\frac{x\log x...\log_{K-2}x(\log_{K-1}x)^B}
{2x\log x...\log_{K-2}x(\log_{K-1}x)^B+x^2(\log x...\log_{K-2}x(\log_{K-1}^x)^B)'}.
\end{align*}
Note that $x^2(\log x...\log_{K-2}x(\log_{K-1}^x)^B)'=o(x\log x...\log_{K-2}x(\log_{K-1}x)^B)$, as $n\rto.$
Thus,
\begin{eqnarray*}
\lim\limits_{n\rto}\frac{\int_{i_0}^n {x\log x...(\log_{K-1}x)^B}dx}{n^2\log n...\log_{K-2}n(\log_{K-1}^n)^B} =\frac12.
\end{eqnarray*}

 If $K=1$ and $B>-1$, $\int_{i_0}^n x^Bdx=\frac1{B+1}(n^{B+1}-i_0^{B+1})\sim\frac1{B+1}n^{B+1};$ if $K=1$ and $B=-1$, $\int_{i_0}^n x^Bdx=\log n-\log i_0\sim\log n;$ otherwise, if  $K=1$ and $B<-1$,
$\int_{i_0}^n x^Bdx=\frac1{B+1}(n^{B+1}-i_0^{B+1})\rightarrow-\frac1{B+1}i_0^{B+1},$ as $n\rto.$

 Note that $x\mapsto x\log x...\log_{K-2}x(\log_{K-1}x)^B$ is monotone in $[k_0,\infty)$,  for some integer $k_0\ge i_0$ depending on $K\text{ and }B.$  Using the same argument as the proofs of \eqref{int2} and \eqref{int3}, we see that for   $K>1$ or  $K=1,\ B\ge -1$,
\begin{align*}\sum\limits_{i=k_0}^n i\log i...\log_{K-2}i(\log_{K-1}i)^B\sim \int_{k_0}^n {x\log x...\log_{K-2}x(\log_{K-1}x)^B}dx
\end{align*}
 as $n\rto$.
Then by Lemma \ref{prhe}(b), we have
\begin{align*}
 \sum\limits_{i=1}^n\rho_1...\rho_i
 \sim\left\{\begin{array}{rcl}
 &cn^2\log n...\log_{K-2}n(\log_{K-1}^n)^B, & K>1,\\
 &c n^{B+1}, & K=1,B>-1,\\
 &c\log n, & K=1, B=-1,
 \end{array}\right.
 \end{align*}
 as $n\rto$.
Moreover, when $K=1,B<-1$, since $\sum\limits_{i=2}^{\infty} i^B$ converges, then $ \sum\limits_{i=1}^{\infty}\rho_1...\rho_i$ converges by Lemma \ref{prhe}(b).

Therefore, by Lemma \ref{prhe}(b), we can see that
\begin{align*}
P&(M=n,D<\infty)
=\frac1{1+\sum_{j=1}^{n-1}\rho_1...\rho_j}\cdot\frac{ \rho_1...\rho_n}{1+\sum_{j=1}^{n}\rho_1...\rho_j}\\
&\sim\left\{\begin{array}{rcl}
&\frac{c n\log n...(\log_{K-1}n)^B}{[n^2\log n...\log_{K-2}n(\log_{K-1}^n)^B]^2}, & K>1, \\
&\frac{cn^B}{(n^{B+1})^2}, & K=1,B>-1,\\
&\frac{\frac cn}{(\log n)^2}, & K=1, B=-1,\\
 & cn^B, & K=1,B<-1,
\end{array}\right.\\
&
=\left\{\begin{array}{rcl}
&\frac{c}{n^3\log n...\log_{K-2}n (\log_{K-1}^n)^B}, & K>1,\\
&\frac{c}{n^{B+2}}, & K=1,B>-1,\\
&\frac{c}{n(\log n)^2}, & K=1, B=-1,\\
 & cn^B, & K=1,B<-1,
\end{array}\right.
\end{align*}
as $n\rto.$ Case 2 is proved and so is Theorem \ref{mxt}.\qed

\vspace{0.1cm}

\noindent{\large{\bf \Large Acknowledgements:}} The programme is carried out during the outbreak of COVID 2019 epidemic in the world. We would like to thank   all people, especially those medical workers, over the world who fight against the virus. Without their sacrifice and dedication, we cannot be so safe. Also, we are very grateful to Prof. W.M. Hong who introduces to us  the Lamperti problem.

\end{document}